\documentclass{article}

\usepackage{amsmath}
\usepackage{amsfonts}
\usepackage{multirow}
\usepackage{graphics}

\begin{document}

\title{A Family of Hybrid Random Number Generators \\ with Adjustable
Quality and Speed}

\author{William K. Cochran\thanks{Two Sigma Investments, 100 Avenue
of the Americas, 16th Floor, New York, NY 10013}
\and
Michael T. Heath\thanks{Department of Computer Science, University
of Illinois at Urbana-Champaign, 201 North Goodwin Avenue, Urbana,
IL 61801}
\and
Kyle W. McKiou\thanks{Department of Computer Science, University
of Illinois at Urbana-Champaign, 201 North Goodwin Avenue, Urbana,
IL 61801}}

\maketitle

\begin{abstract}
Conventional random number generators provide the speed but not
necessarily the high quality output streams needed for large-scale
stochastic simulations.  Cryptographically-based generators, on the
other hand, provide superior quality output but are often deemed
too slow to be practical for use in large simulations.  We combine
these two approaches to construct a family of hybrid generators
that permit users to choose the desired tradeoff between quality
and speed for a given application.  We demonstrate the effectiveness,
performance, and practicality of this approach using a standard
battery of tests, which show that high quality streams of random
numbers can be obtained at a cost comparable to that of fast
conventional generators.
\end{abstract}



\section{Random Number Generators}

Large stochastic simulations, such as those typically run on highly
parallel supercomputers, consume vast quantities of random numbers.
Designers and users of random number generators face a dilemma:
how to achieve the speed necessary to make the computation feasible
while avoiding any hint of correlation or bias that would invalidate
the results.  Conventional random number
generators~\cite{anderson90,james90,lecuyer90}, typically based on
simple arithmetic recurrences, may possess the requisite speed, but
they carry a significant risk of failure due to the questionable
quality of the resulting streams, especially when produced in great
abundance, thereby inviting more opportunity for detectable patterns
to emerge~\cite{hellekalek98,marsaglia68,park_miller88}.  This should
not be surprising: according to Kolmogorov's definition~\cite{kolmogorov65},
a truly random stream is \emph{incompressible} (i.e., it has no shorter
description than simply enumerating its members); a generator
algorithm provides a succinct encoding of the stream it produces,
which therefore cannot be truly random.  Such deterministically
generated (hence reproducible) streams are more properly called
\emph{pseudorandom}, but nevertheless a well-designed generator may
still produce streams that appear random in that they may pass
various statistical tests of randomness~\cite{lecuyer07,marsaglia02}.

\subsection{Cryptographically-Based Generators}

An alternative view of randomness is taken in cryptography, in which
random number generators are ubiquitously employed to generate
random keys, signatures, and the like~\cite{goldreich01}.
In security applications, quality is paramount and performance
secondary, as an easily compromised cryptographic protocol would
be useless regardless of its efficiency.
Here, a stream is considered random if it is
\emph{computationally indistinguishable} from the desired true
distribution, that is, detecting any difference would require a
prohibitive amount of computation~\cite{goldreich10}.  Many
cryptographic protocols are based on \emph{one-way functions}, that
is, functions that are easy to compute but difficult to invert (i.e,
to find an input value that yields a given output value).  In
principle, a one-way function can be used to construct a ``perfect''
random number generator, in that resistance to inversion ensures
computational indistinguishability of the resulting
stream~\cite{goldreich10,hastad1999}.  Alas, there are two difficulties
with this putative panacea: the existence of one-way functions is
an open theoretical question (closely related to the famously
unresolved $P = \mathit{NP}$?  problem)~\cite{levin03}, and the
known candidates believed to be one-way functions (e.g., factoring
products of large primes) are relatively expensive to evaluate,
making the resulting generator slow and thus seemingly unsuitable
for large stochastic simulations~\cite{lecuyer_proulx89}.  The first
of these issues does not seem to be a problem in practice unless
efficient algorithms are unexpectedly discovered for these currently
intractable inversion problems, or if it turns out that
$P = \mathit{NP}$, in which case true one-way functions do not
exist.  The speed of cryptographically-based generators is an
important practical issue, however, which we address next.

Practical cryptographic transforms, such as block ciphers and hash
functions, are based on algorithms that are carefully designed and
extensively tested to ensure their resistance to inversion.  These
algorithms employ basic operations (e.g., substitutions, permutations)
that obscure any structure in the input by dissipating it over the
output.  This process is usually repeated multiple (often many)
times (called ``rounds'') to ensure sufficient diffusion, which is
the main reason such transforms are relatively slow.  Nevertheless,
their superior ability to produce high quality random output is
crucial in security applications, and it could also make stochastic
simulations more robust if the speed of cryptographically-based
generators could be made competitive with conventional generators.
Recently, cryptographically-based generators have been used for a
few stochastic simulations and computer graphics applications.
Cryptographic transforms employed have included AES and
Threefish~\cite{salmon_etal11}, MD5~\cite{tzeng_wei08}, and
TEA~\cite{zafar_olano_curtis10}.  To make the cryptographically-based
generators competitive in speed, the corresponding cryptographic
transforms were weakened by reducing the number of rounds significantly,
and in some cases by making other simplifications to the algorithms
as well.  The resulting generators are no longer sufficiently secure
for cryptographic applications, but still produce streams of
sufficient quality for stochastic simulations, as demonstrated by
passing standard statistical test batteries.  Here we will pursue
the same objective but with a different approach, using a full-strength
cryptographic transform but combining the output of the resulting
generator with that of a conventional generator, which will allow
us explicit, quantitative control of the tradeoff between quality
and speed.

\section{Combined Cryptographic/LCG Generator}

To create a fast generator with sufficient randomness, we combine
a crypto\-graphically-based generator with a fast conventional
generator.  Because of their speed, simplicity, and widespread use,
we chose a linear congruential generator (LCG) for the conventional
generator (which will also serve to illustrate that even a relatively
low-quality generator can be rehabilitated using our technique),
but any type of conventional generator should also work.  Because
cryptographic transforms are relatively slow, we amortize the cost
of computing random bits by reusing them in the combined generator.
We define two parameters relevant to combining the two generators:
\begin{itemize}
\item
\emph{size}: the number of cryptographically-based 32-bit random
numbers (cryptographic values) stored at a given time;
\item
\emph{repetition}: the number of times each cryptographic value is
used before being discarded.
\end{itemize}

\begin{figure}
{\resizebox*{4.8in}{!}{\includegraphics{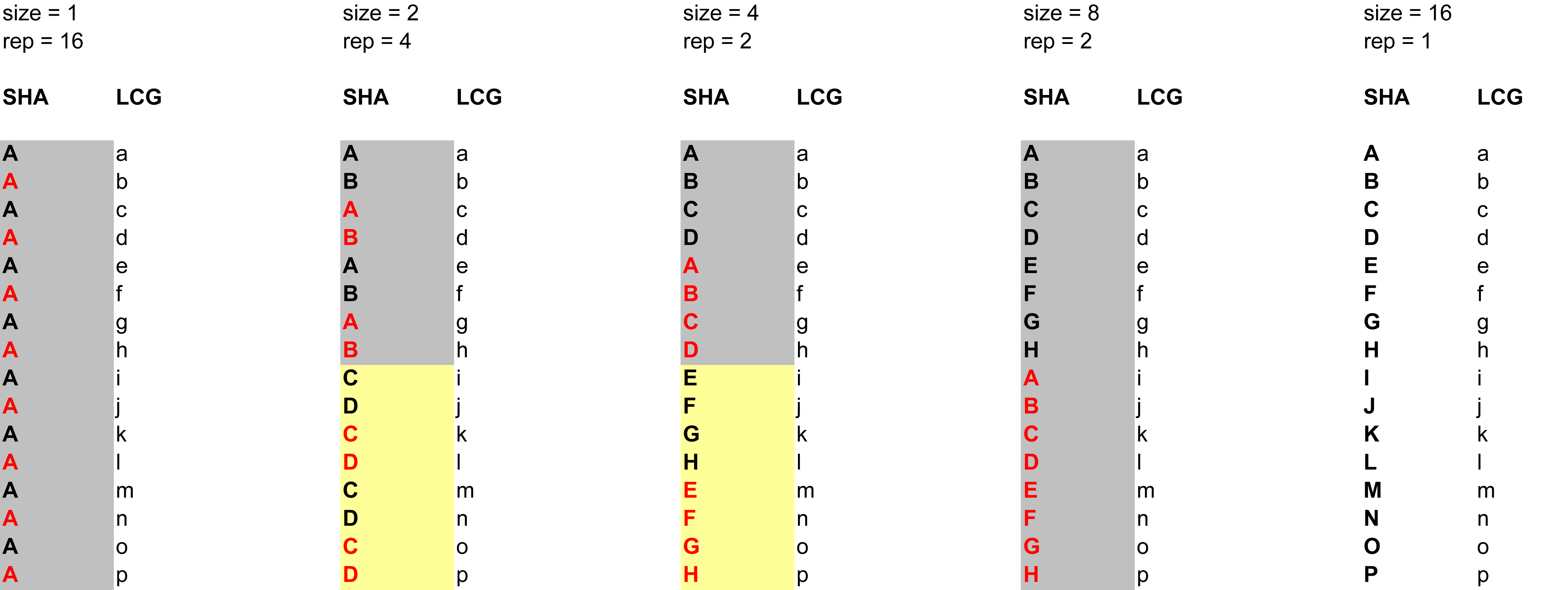}} \par}
\caption{Examples showing effects of size and repetition parameters in
combined generator.}
\label{fig1}
\end{figure}

Motivated by its well-known bias-reducing property~\cite{davies02},
we use bitwise exclu\-sive-or (XOR) to combine cryptographic
values with LCG-based 32-bit random numbers (LCG values) in an
alternating fashion.  Let $k$ be the value of size and $n$ the value
of repetition.  The first $k$ LCG values are XORed uniquely with
the $k$ cryptographic values.  Specifically, the first LCG value
is XORed with the first cryptographic value, the second LCG value
is XORed with the second cryptographic value, and so on, until the
$k$th LCG value is XORed with the $k$th cryptographic value.  The
next $k$ LCG values are XORed with the same $k$ cryptographic values
in a similar fashion.  This process is repeated until each cryptographic
value has been used to XOR $n$ LCG values. Then the $k$ cryptographic
values are discarded and $k$ new cryptographic values are produced.
This procedure is repeated until the generator stream is terminated.
Figure~\ref{fig1} illustrates the effects of various values of the
size and repetition parameters in the combined generator for some
small examples.

\subsection{Practical Implementation}

We chose the Secure Hash Algorithm SHA-256~\cite{fips180-4} as the
basis for the cryptographic generator due to its exceptional
resistance to inversion and its wide acceptance in the cryptographic
community.  Because each call to SHA-256 produces 256 random bits
of output, eight 32-bit cryptographic values are produced by each
call.  Due to the cryptographic strength of SHA-256, we use
$i = 1, 2, 3, \ldots$ to seed the generator (as in~\cite{salmon_etal11}),
with no noticeable correlation in the resulting output.

We selected three popular LCGs to combine with SHA-256:
\begin{itemize}
\item
Super-Duper~\cite{marsaglia72}:
$x_n = x_{n-1}*69069 + 1 \pmod{2^{32}}$,
\item
glibc \texttt{rand}:
$x_n = x_{n-1}*1103515245+ 12345 \pmod{2^{32}}$,
\item
Borland C++ \texttt{rand}:
$x_n = x_{n-1}* 22695477 + 1 \pmod{2^{32}}$.
\end{itemize}
Because the low-order bits in power-of-two modulus LCGs have much
shorter periods than the high-order bits~\cite{mascagni00},
we call each LCG twice and concatenate the 16 highest-order
bits of each iteration to obtain a more reliable 32-bit random
number.

\subsection{Parallel Generators}

One of the key advantages of using a cryptographic transform as
part of our combined generator is ease of parallelization.  SHA-256
has an exceptionally large $(2^{64}-1)$-bit seed space~\cite{fips180-4},
which allows for production of $8*2^{2^{64}-1}/n$ independent streams
of $n$ 32-bit random numbers.  For example, the SHA-256 generator
can easily produce $2^{1000000}$ independent streams of $2^{1000000}$
32-bit random values if desired.
To parallelize the SHA-256 generator, simply concatenate the
stream number with a counter $i=1, 2, 3, \ldots$ for each stream.
Because each stream will have a unique seed, the cryptographic
strength of SHA-256 ensures independence of streams.  To parallelize
the LCG, one can choose any of several suitable ways as outlined
by Srinivasan, Mascagni, and Ceperley~\cite{srinivasan03}, such as
introducing a lag, choosing prime addends, or choosing prime moduli.

\section{Test Results}

To test our combined generators, we used the test batteries in
TestU01~\cite{lecuyer07}, a collection of popular random number
generators and tests to measure their effectiveness.  TestU01
includes three predefined test batteries for random bit sequences:
SmallCrush, Crush, and BigCrush.  We used SmallCrush for preliminary
testing, but it is insufficiently stringent to detect anything
other than gross defects in generators.  Crush is a considerably
more strenuous battery of tests, applying 96 statistical tests and
producing 144 test statistics and $p$-values while consuming
approximately $2^{35}$ random numbers.  Crush was applied to all
generators in this paper to assess their approximate quality and
suggest candidates for further testing.  BigCrush was then used to
test the most promising generators more rigorously.  BigCrush applies
106 statistical tests and produces 160 test statistics and $p$-values
while consuming approximately $2^{38}$ random numbers.  All programs
were written in C, compiled with gcc, and run on a 2.8 GHz Intel
Xeon X5560 processor.

\begin{figure}
{\resizebox*{4.8in}{!}{\includegraphics{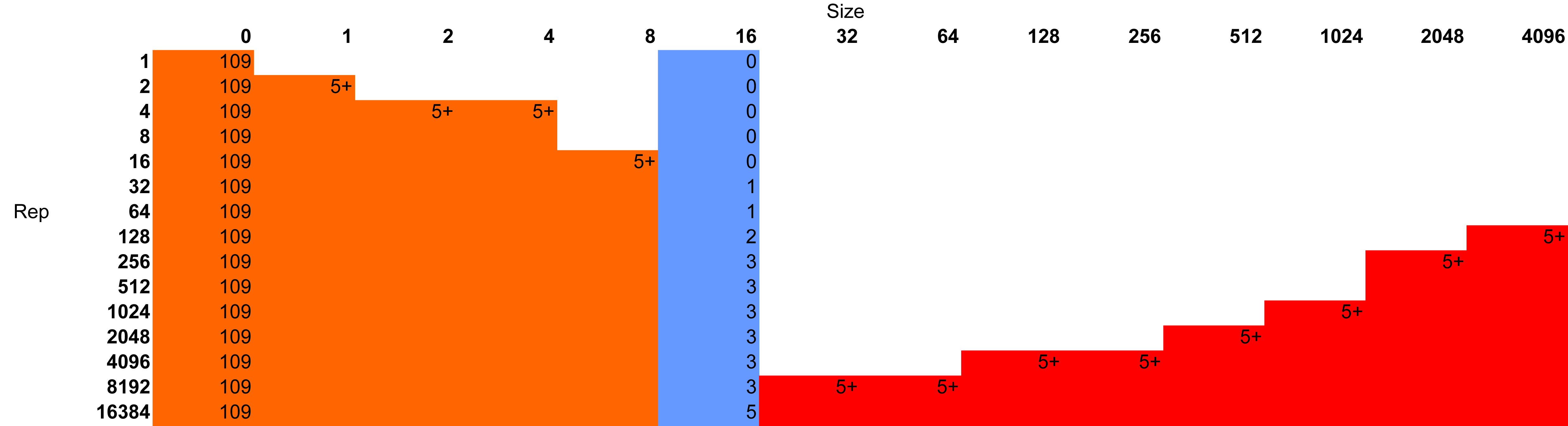}} \par}
\caption{BigCrush failures (out of 160 total) for combined SHA-256/Super-Duper
generator.}
\label{fig2}
\end{figure}

\begin{figure}[b!]
\begin{minipage}[t]{2.4in}
{\centering \resizebox*{2.4in}{!}{\includegraphics{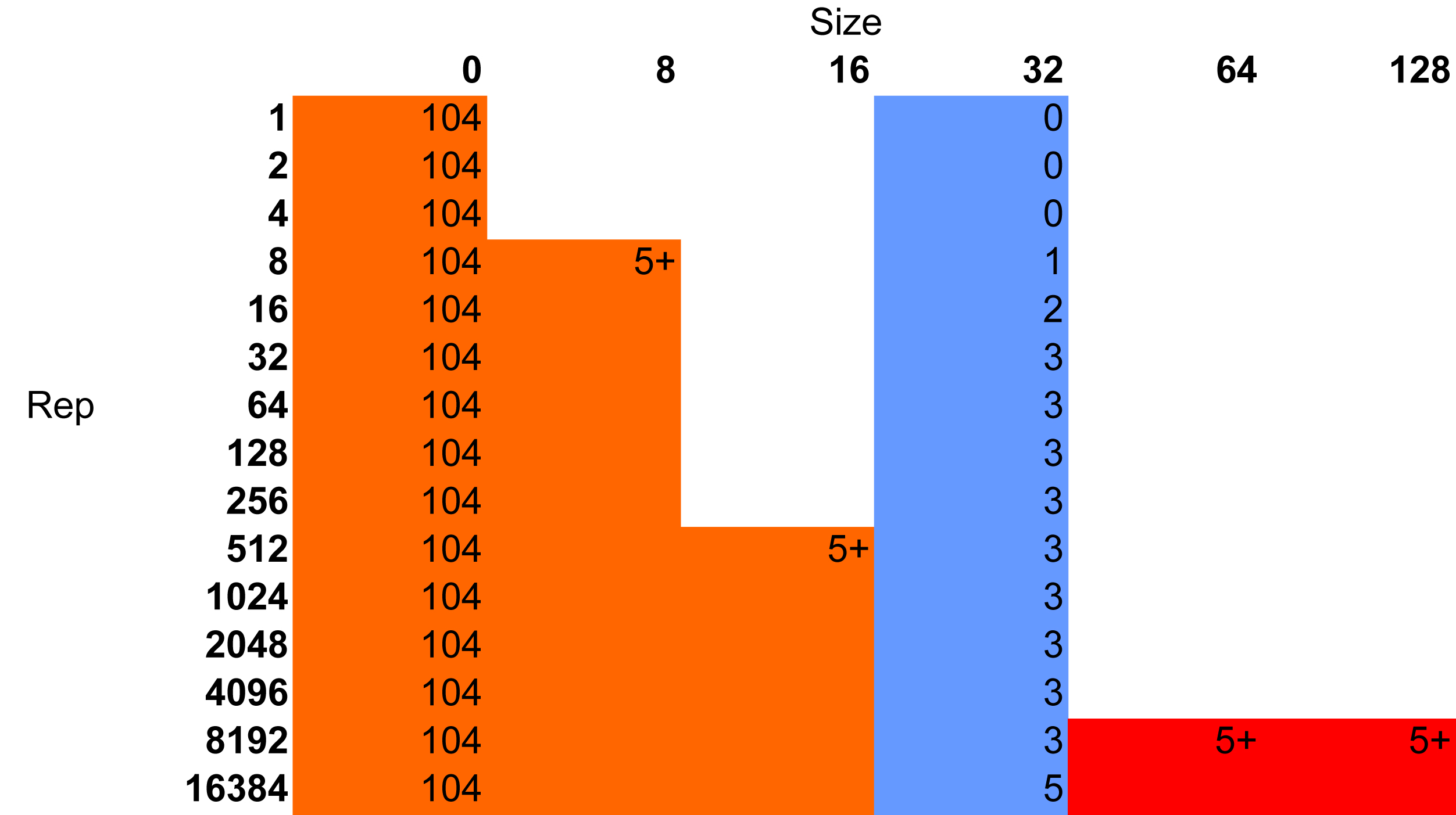}} \par}
\end{minipage}
\begin{minipage}[t]{2.4in}
{\centering \resizebox*{2.12in}{!}{\includegraphics{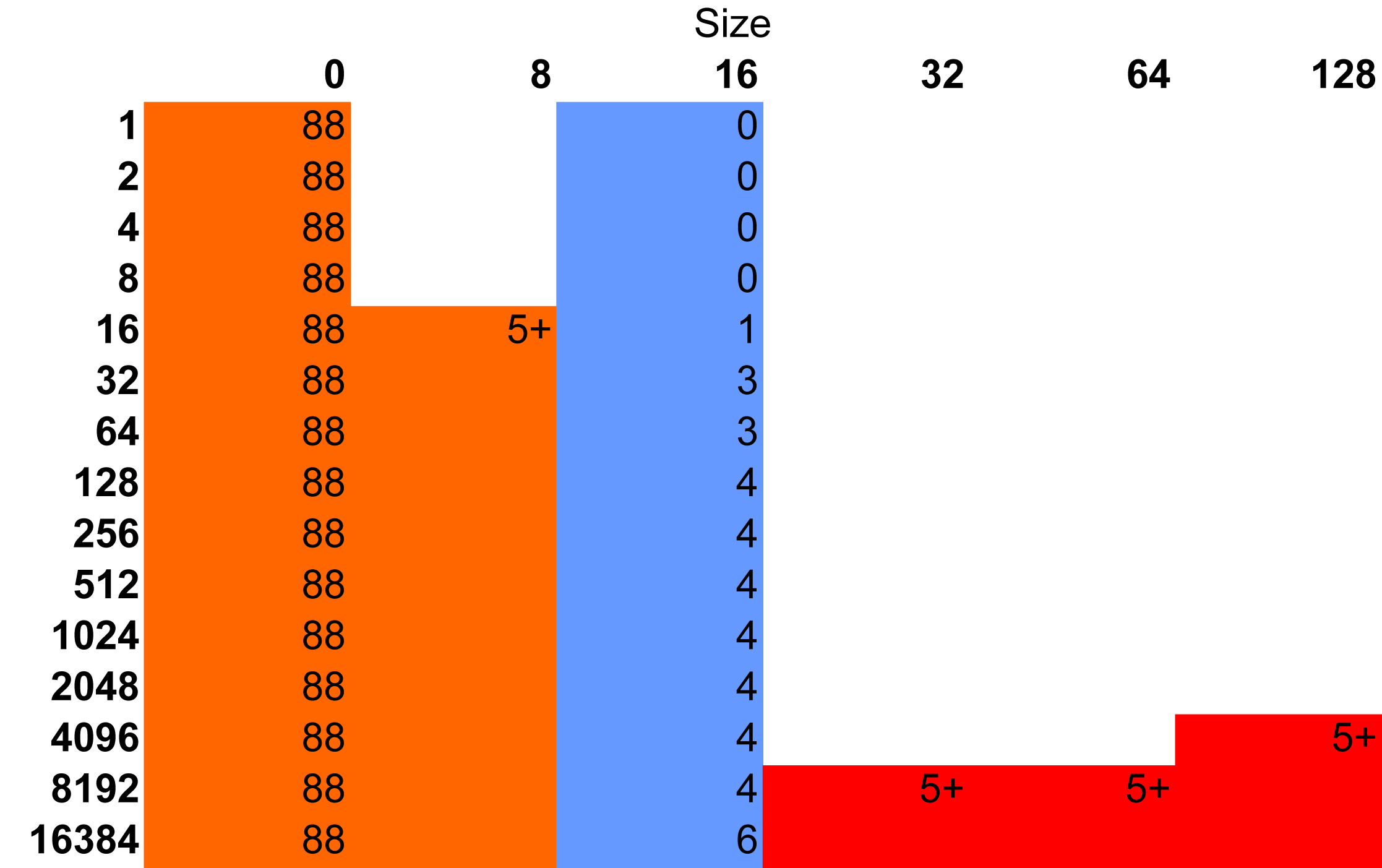}} \par}
\end{minipage}
\caption{BigCrush failures (out of 160 total) for combined SHA-256/glibc
\texttt{rand} generator (left) and combined SHA-256/Borland C++ \texttt{rand}
generator (right).}
\label{fig3}
\end{figure}

Figures~\ref{fig2} and \ref{fig3} show results of BigCrush tests
for each of the three combined generators for various values of the
size and repetition parameters.  The column for size=0 gives the
number of test failures (out of 160) using only the corresponding
LCG, for which repetition is irrelevant.  For a given nonzero value
of size, increasing repetition (reuse of cryptographic values)
degrades reliability (larger numbers of test failures), as expected.
Varying the size parameter reveals a more complex relationship, in
which increasing size from $1$ improves reliability up to an optimal
value, then may degrade reliability beyond that point.  An explanation
for this behavior is given in the discussion below.  To avoid running
expensive BigCrush tests for cases that were certain to have a
relatively large number of failures (and therefore not yield a
useful generator), the figures merely indicate (by orange and red
shading) repetition factors for which the number of failures is
$5$ or more (out of 160), except that full results are given for
the optimal value of size for each generator (blue shading).  The
three generators produced generally similar results, although the
SHA-256/glibc \texttt{rand} combined generator achieves maximum
reliability for size=32, whereas size=16 is optimal for the other
two combined generators.

Timing results given in Table~\ref{timings} show that increasing
repetition increases the speed of the combined generator, i.e., the
more times the cryptographic values are reused, the faster the
generator produces random bits.  In Table~\ref{timings}, we see
that the SHA-256 generator takes an average of 1.56 seconds to
produce $2^{23}$ random values, whereas each of the LCGs alone takes
an average of 0.19 seconds to produce the same number of values.
The combined generator that uses each cryptographic value only once
(repetition=1) slows performance due to the overhead of creating
both cryptographic values and LCG values.  However, the speed is
improved by producing fewer cryptographic values and increasing the
number of times they are used in the combined generator stream.  As
repetition increases, speed initially increases substantially, but
improvement levels off for repetition $>$ 64.  The combined generators
with repetition $\geq 256$ are nearly identical in speed to the LCG
generator.  Table~\ref{timings} gives timings only for the optimal
value of size for each combined generator, but additional tests
showed that for a given fixed repetition, timing results remain
nearly identical as the size parameter varies.

\begin{table}[tbh]
\centering
\caption{Execution times in seconds (averaged over 20 trials) for producing
$2^{23}$ 32-bit random numbers for three combined generators, each
combining SHA-256 with indicated LCG.  Failures indicates number of test
statistics failed out of 160 using BigCrush.}
\label{timings}
\vspace{0.5cm}
\begin{tabular}{r | r r | r r | r r} \hline
\multicolumn{1}{c}{LCG} &
\multicolumn{2}{c}{Super-Duper} &
\multicolumn{2}{c}{glibc \texttt{rand}} &
\multicolumn{2}{c}{Borland C++ \texttt{rand}} \\
\multicolumn{1}{c}{size} &
\multicolumn{2}{c}{16} &
\multicolumn{2}{c}{32} &
\multicolumn{2}{c}{16} \\ \hline
\multicolumn{1}{c}{repetition} &
\multicolumn{1}{c}{time} &
\multicolumn{1}{r}{failures} &
\multicolumn{1}{c}{time} &
\multicolumn{1}{r}{failures} &
\multicolumn{1}{c}{time} &
\multicolumn{1}{r}{failures}
\\ \hline
SHA & 1.56 & 0 & 1.56 & 0 & 1.56 & 0 \\
1 & 1.70 & 0 & 1.70 & 0 & 1.70 & 0 \\
2 & 0.95 & 0 & 0.95 & 0 & 0.94 & 0 \\
4 & 0.57 & 0 & 0.57 & 0 & 0.57 & 0 \\
8 & 0.38 & 0 & 0.38 & 1 & 0.38 & 0 \\
16 & 0.28 & 0 & 0.28 & 2 & 0.28 & 1 \\
32 & 0.24 & 1 & 0.24 & 3 & 0.24 & 3 \\
64 & 0.21 & 1 & 0.21 & 3 & 0.21 & 3 \\
128 & 0.20 & 2 & 0.20 & 3 & 0.20 & 4 \\
256 & 0.19 & 3 & 0.19 & 3 & 0.19 & 4 \\
512 & 0.19 & 3 & 0.19 & 3 & 0.19 & 4 \\
1024 & 0.19 & 3 & 0.19 & 3 & 0.19 & 4 \\
2048 & 0.19 & 3 & 0.19 & 3 & 0.19 & 4 \\
4096 & 0.19 & 3 & 0.19 & 3 & 0.19 & 4 \\
8192 & 0.19 & 3 & 0.19 & 3 & 0.19 & 4 \\
16384 & 0.19 & 5 & 0.19 & 5 & 0.19 & 6 \\
LCG & 0.19 & 109 & 0.19 & 104 & 0.19 & 88 \\ \hline
\end{tabular}
\end{table}

\subsection{Discussion}

Reusing SHA-256 values in the combined generators results in less
costly but less reliable streams of random numbers.  Table~\ref{timings}
illustrates the relationship between repetition and both speed and
failure rate for each of the combined generators.  For example, the
combined SHA-256/Super-Duper generator with repetition=16 has
a more than five-fold increase in speed compared to the SHA-256
generator alone, while still passing \emph{all} of the statistical
tests.  At repetition=256, the combined generator is essentially
identical in speed to the LCG generator, while failing only 3/160
tests compared to failing 109/160 tests using Super-Duper alone.

Varying size produces slightly less intuitive results.  It is well
known that LCGs exhibit a high degree of serial
correlation~\cite{hellekalek98,marsaglia68,park_miller88}.  With a
relatively small size parameter, we are able to eliminate correlation
between successive LCG values.  Unfortunately, correlation remains
in each substream that shares a common cryptographic-value buffer.
For example, if we have a high-repetition generator with size=2,
there will be significant correlation between every second value
in the sequence.  This is the same type of correlation that lagged
LCGs suffer from.  Because this correlation stems from the fact
that each LCG value is determined uniquely by the previous value
(given that the other parameters are fixed), we call this type
of correlation \emph{local} correlation (indicated by orange shading
in Figures~\ref{fig2}--\ref{fig3}).

Increasing size improves generator quality, but only up to a point.
Clearly, as size increases, local correlation in the resulting
stream decreases.  But as size becomes large, another type of
correlation in the LCG stream becomes visible between long, disjoint
sequences of values.  This correlation is a result of the simple
recursive nature of the LCG, leading to an inherent structure in
the stream regardless of the actual values produced~\cite{mascagni00}.
We call this type of correlation \emph{global} correlation due to
its ubiquitous nature (indicated by red shading in
Figures~\ref{fig2}--\ref{fig3}).

Thus, the size parameter should be chosen sufficiently large to
hide local correlation, but not so large as to reveal global
correlation in the LCG stream.  Our results suggest that combined
generators with size between 16 and 32 are most suitable in this
regard (blue shading in Figures~\ref{fig2}--\ref{fig3}).

\section{Recommendations}

Several factors should be considered when choosing a combined random
number generator.  We will start with the easiest parameter to
choose: size.  Clearly, it is preferable to choose the size that
minimizes failure rate.  If several size values produce similarly
reliable generators, then the smallest acceptable size should be
chosen to minimize the amount of state required, which in turn
minimizes the amount of memory required.  This is particularly
important for highly-parallel applications, for which careful use
of memory is often critical.  If minimizing state is absolutely
essential, a suboptimal size can be chosen at some cost in reliability.

To choose the best repetition factor, users must carefully consider
their needs.  Consider, for example, the SHA-256/Super-Duper combined
generator.  In Table~\ref{timings}, we see that the repetition $\leq
16$ generators provide the highest reliability, while the repetition
$\geq 256$ generators provide maximum speed.  Choosing the repetition
parameter may be easy for users who are either very conservative
($\leq 16$) or care mainly about speed ($\geq 256$), but other users
may prefer a compromise between speed and reliability.  For example,
many users may find the size=128 generator suitable for their needs,
as it passes all but 2/160 BigCrush tests, and thus is vastly
superior in quality to standard library LCGs, yet runs only 3.2\%
slower than the pure LCG.  Because it is easy to vary the repetition
parameter, users can experiment and choose the value that best meets
their needs in a given situation, for example, ``quick and dirty''
preliminary exploration versus more exacting final simulation runs.

When choosing a cryptographic transform as the basis for the
cryptographic generator, we recommend SHA-256 for reasons already
given, although other cryptographic functions may be equally suitable,
and quite possibly faster, which might give somewhat different
tradeoffs than we observed in our tests.  For relatively reliable
and easily parallelizable LCGs, we recommend the SPRNG
library~\cite{mascagni00}.  These generators have been extensively
studied and are suitable for highly parallel implementation, and
based on our results they should perform well as part of a size=16
or size=32 combined generator.

\section{Conclusion}

In this paper we have shown that the use of a full-strength
cryptographic transform in a combined generator can dramatically
improve the quality of streams produced by a conventional linear
congruential generator for stochastic simulations, with no significant
sacrifice in speed.  Moreover, our approach provides explicit,
quantitative control of the tradeoff between speed and reliability,
so that users can easily select the choice that best meets their needs
in a particular situation.  Future work along this line could include
consideration of other cryptographic transforms besides SHA-256,
as well as other types of conventional generators besides LCGs,
such as lagged Fibonacci generators, for the respective components
of a combined generator.  Either or both of these options may produce
somewhat different tradeoffs between speed and reliability from those
we observed in our experiments.

\bibliographystyle{plain}
\bibliography{refs}

\end{document}